\documentclass[a4paper,10pt]{amsart}

\usepackage[T1]{fontenc} \usepackage[utf8]{inputenc}
\usepackage{lmodern}

\usepackage{graphicx}

\usepackage{enumerate}

\usepackage{amsfonts,amsmath,amstext,amsbsy,amssymb}
\usepackage{amsbsy,amsopn,amsthm,amscd,amsxtra}

\usepackage{latexsym,mathrsfs}

\usepackage{mathabx} \newcommand{\carr}{\righttoleftarrow}

\usepackage{hyperref}
\RequirePackage[dvipsnames]{xcolor} % [dvipsnames]
\definecolor{halfgray}{gray}{0.55} 
\definecolor{webgreen}{rgb}{0,0.5,0}
\definecolor{webbrown}{rgb}{.6,0,0} \hypersetup{%
  colorlinks=true, linktocpage=true, pdfstartpage=3,
  pdfstartview=FitV,%
  breaklinks=true, pdfpagemode=UseNone, pageanchor=true,
  pdfpagemode=UseOutlines,%
  plainpages=false, bookmarksnumbered, bookmarksopen=true,
  bookmarksopenlevel=1,%
  hypertexnames=true,
  pdfhighlight=/O,%hyperfootnotes=true,%nesting=true,%frenchlinks,%
  urlcolor=webbrown, linkcolor=RoyalBlue,
  citecolor=webgreen, %pagecolor=RoyalBlue,%
  pdftitle={},%
  pdfauthor={},%
  pdfsubject={2000 MAthematical Subject Classification: Primary:},%
  pdfkeywords={},%
  pdfcreator={pdfLaTeX},%
  pdfproducer={LaTeX with hyperref}%
}

\usepackage[T1]{fontenc} \usepackage{ae}

\newcommand{\abs}[1]{\left\lvert{#1}\right\rvert}
\newcommand{\norm}[1]{\left\|{#1}\right\|}

\newcommand{\bb}{\mathbb}

\newcommand{\R}{\mathbb{R}}\newcommand{\N}{\mathbb{N}}
\newcommand{\Z}{\mathbb{Z}}
\newcommand{\T}{\mathbb{T}}\newcommand{\C}{\mathbb{C}}

 \newcommand{\ie}{i.e.\ }
\newcommand{\etc}{\emph{etc.\ }}

\newtheorem{theorem}{Theorem}[section]
\newtheorem{proposition}[theorem]{Proposition}
\newtheorem{lemma}[theorem]{Lemma}
\newtheorem{corollary}[theorem]{Corollary}

\theoremstyle{definition} \newtheorem{definition}{Definition}

\newtheorem{question}[theorem]{Question}

\theoremstyle{remark} \newtheorem{remark}[theorem]{Remark}

\newcommand{\Diff}[2]{\mathrm{Diff}^{#1}_{#2}}
\newcommand{\Homeo}[1]{\mathrm{Homeo}_{#1}}
\newcommand{\Dis}{\mathcal{D}} \newcommand{\Deck}{\mathrm{Deck}}
\newcommand{\Isom}{\mathrm{Isom}}

 \newcommand{\dd}{\:\mathrm{d}}

\DeclareMathOperator{\SL}{SL} \DeclareMathOperator{\SO}{SO}
\DeclareMathOperator{\GL}{GL} 
\DeclareMathOperator{\PSL}{PSL}

 \newcommand{\Ss}{\mathbb{S}}
 \newcommand{\e}{\varepsilon}

\title[The Burnside problem]{The Burnside problem for
  $\Diff{}\omega(\bb{S}^2)$}

\author[Hurtado]{Sebastián Hurtado} \address{University of Chicago,
  Illinois, USA}
\email{shurtado@chicago.edu}

\author[Kocsard]{Alejandro Kocsard} \address{Universidade Federal
  Fluminense, Rio de Janeiro,
  Brasil}\urladdr{www.professores.uff.br/kocsard}
\email{akocsard@id.uff.br}

\author[Rodríguez-Hertz]{Federico Rodríguez-Hertz}
\address{Pennsylvania State University, USA}
\email{hertz@math.psu.edu}

\date{\today}

\begin{document}

\begin{abstract}
  A group $G$ is periodic of bounded exponent if there exists
  $k \in \mathbb{N}$ such that every element of $G$ has order at most
  $k$. We show that every finitely generated periodic group of bounded
  exponent $G<\Diff{}{\omega}(\bb{S}^2)$ is finite, where
  $\Diff{}{\omega}(\bb{S}^2)$ denotes the group of diffeomorphisms of
  $\bb{S}^2$ that preserve an area-form $\omega$.
\end{abstract}

\maketitle

\section{Introduction}
\label{sec:intro}

A group $G$ is said to be \emph{periodic} if every element of $G$ has
finite order. If there exists $N \in \mathbb{N}$ such that $g^N = id$
for every $g \in G$ (where $id\in G$ denotes the identity element),
then $G$ is said to be a periodic group of \emph{bounded
  exponent}. The so called \emph{Burnside problem} is a famous
question in group theory originally considered by Burnside in
\cite{BurnsideUnsettledQuestion} which can be stated as follows:

\begin{question}[Burnside, 1905]
  \label{ques:Burnside-general}
  Let $G$ be a finitely generated periodic group. Is $G$ necessarily
  finite? And if $G$ is periodic of bounded exponent?
\end{question}

Burnside himself proved in his article
\cite{BurnsideUnsettledQuestion} that if $G$ is a linear finitely
generated periodic group of bounded exponent (\ie $G<\GL_n(\C)$), then
$G$ must be finite.  In 1911, Schur~\cite{SchurGruppenLinearer}
improved Burnside's result removing the bounded exponent hypothesis.

However, in general the answer to Question~\ref{ques:Burnside-general}
turned out to be ``no'' as counterexamples were later discovered by
Golod-Schafarevich\cite{GolodNilAlg, Golod1ShafaClassField} and
Adian-Novikov\cite{NovikovAdjanInfinitePeriodic} in the 1960's. Since
then, many more examples have been constructed by Olshanskii, Ivanov,
Grigorchuk among others, and there is a vast literature on the subject
(see for instance \cite{BurnsideWebAnon}).

``Non-linear'' transformation groups like homeomorphism groups,
diffeomorphisms groups, volume preserving diffeomorphism groups,
groups of symplectomorphisms, \etc are conjectured to have many common
features with linear groups (see Fisher's survey
\cite{FisherRecentZimmerProg} on the Zimmer program). For example, the
following question is attributed to E. Ghys and B. Farb is stated in
\cite{FisherRecentZimmerProg} (see also \cite[Question
13.2]{FisherGroupsActing}):

\begin{question}[Burnside problem for homeomorphism groups]
  \label{ques:Burnside-homeo-groups}
  Let $M$ be a connected compact manifold and $G$ be a finitely
  generated subgroup $G<\Homeo{}(M)$ such that all the elements of $G$
  have finite order. Then, is $G$ necessarily finite?
\end{question}

At this point it is important to remark that compactness is an essential
hypothesis in Question~\ref{ques:Burnside-homeo-groups}. In fact, it
is well known that given any finitely presented group $G'$, there
exists a connected smooth manifold $M$ such that its fundamental group
$\pi_1(M)$ is isomorphic to $G'$. Hence, (any sub-group of) $G'$
clearly acts faithfully on the universal cover $\tilde M$ of $M$. To
the best of our knowledge, so far it is not known whether there exists
a finitely presented infinite periodic group. However,  free burnside groups and Grigorchuk's group \cite{grigorchuk1980burnside} are known to be recursively presented, and consequently, can be embedded by Higman's embedding Theorem into a finitely presented group $G'$.

For the time being, we know Question~\ref{ques:Burnside-homeo-groups}
has a positive answer just in a few cases and no negative one is
known.

For instance, in the one-dimensional case, \ie when $M=\Ss^1$, this is
an easy consequence of Hölder's theorem (see for instance
\cite[Theorem 2.2.32]{NavasGroupsCircleDiff}): in fact, one can show
that any periodic group $G<\Homeo+(\Ss^1)$ acts freely on $\Ss^1$ (\ie
the identity is the only element exhibiting fixed points), and by
Hölder's theorem, $G$ is conjugate to a subgroup of $\SO(2)$. In
particular, $G$ is abelian and so, it is finite if and only if is
finitely generated.

In higher dimensions, Rebelo and Silva~\cite{RebeloBurnside} give a
positive answer to Question~\ref{ques:Burnside-homeo-groups} for
groups of symplectomorphisms on certain symplectic $4$-manifolds; and
Guelman and Lioussse \cite{GuelmanLiousseBurnsideMeasurePres,
  GuelmanLiousseBurnsideGroupsHomeo} for hyperbolic surfaces and
groups of homeomorphisms of $\T^2$ exhibiting an invariant probability
measure. In \S\ref{sec:burnside-problem-Td-hyperbolic-mfds} we extend
these last results for higher dimensions. 

The main results of this article are the following:

\begin{theorem}
  \label{thm:main}
  Let $\omega$ be an area form on $\Ss^2$ and
  $\Diff\infty\omega(\bb{S}^2)$ be the group
  $C^\infty$-diffeomorphisms of $\Ss^2$ that preserve $\omega$.

  Then, any finitely generated periodic subgroup of bounded exponent of
  $\Diff\infty\omega(\Ss^2)$ is finite.
\end{theorem}

\begin{theorem}
  \label{thm:Burnside-problem-hyperbolic-mfd}
  Let $M^q$ be a compact hyperbolic manifold of dimension $q\geq
  3$. Then, any periodic subgroup of $\Homeo{}(M)$ is finite.
\end{theorem}

Observe that in Theorem~\ref{thm:Burnside-problem-hyperbolic-mfd} we
are not \emph{à priori} assuming the subgroup is finitely generated.

\begin{theorem}
  \label{thm:Burnside-problem-Td}
  Let $\mu$ be a Borel probability measure on the $q$-dimensional
  torus $\T^q$, and $\Homeo\mu(\T^q)$ denote the group of
  homeomorphisms of $\T^q$ that preserves the measure $\mu$.

  Then, any finitely generated periodic subgroup of $\Homeo\mu(\T^q)$
  is finite.
\end{theorem}

We should point out that Conejeros\cite{conejeros2018periodic} recently proved some results about periodic groups in $\Homeo+(\Ss^2)$, using the theory of rotation sets and by methods different than ours and which are closer to the methods used by Guelman and Lioussse \cite{GuelmanLiousseBurnsideMeasurePres,
  GuelmanLiousseBurnsideGroupsHomeo}.

\subsection{Outline of the paper}
\label{sec:outline}

In \S\ref{sec:preliminaries}, we fix some notations we shall use all
along the paper and recall some previous known results on differential
geometry and topology.

Sections \ref{sec:derivatives}, \ref{sec:sub-expon-growth-case} and
\ref{sec:exponential-growth-case} are dedicated to the proof of
Theorem~\ref{thm:main}. For the sake of readability, let us roughly
explain the strategy of the proof.

To do that, let $G$ be a finitely generated periodic sub-group of
$\Diff\infty\omega(\Ss^2)$ with bounded exponent. In
\S\ref{sec:derivatives} we show that the group $G$ exhibits
\emph{sub-exponential growth of derivatives,} \ie the norm of the
derivatives of the elements of $G$ grows sub-exponentially with
respect to their word length in $G$. The proof of the sub-exponential
growth of derivatives is based on encoding the group action as a
$\Diff{}{}(\Ss^2)$-cocycle over the full-shift $k$ symbols, where $k$
denotes the amount of generators of $G$. This is a classical
construction in random dynamics (see for instance
\cite{kifer2006random}) and we use it to show that the
exponential growth of derivative implies the existence of an
element of $G$ exhibiting a hyperbolic periodic point, which
contradicts the fact that $G$ is periodic. These ideas are closely
related Livsic's theorem for diffeomorphism group cocycles
\cite{KocsardPotrieLivThm, AvilaKocLiu} and to Katok's closing lemma
(see Lemma~\ref{Katok}).

Then, invoking a rather elementary argument, in
\S\ref{sec:exponential-growth-case} we show that the sub-exponential
growth of derivatives is incompatible with the exponential growth of
the group $G$ itself, with respect to the word length. It is
interesting to remark that this is the only part of the proof where
boundedness of the exponent is indeed used.

So, in \S\ref{sec:sub-expon-growth-case}, which contains the more
elaborate and intricate arguments of the proof, we assume the group
$G$ has sub-exponential growth (with respect to the word length).
Under this hypothesis, we show the existence of a diffeomorphism
$h\colon\Ss^2\carr$ such that $h^{-1}Gh<\SO(3)$. Then, as consequence
Schur's theorem \cite{SchurGruppenLinearer}, $G$ must be finite.

The existence of such a diffeomorphism $h$ will be given by the
Uniformization theorem of surfaces
(Theorem~\ref{thm:uniformization}). In fact, we will construct a
$G$-invariant smooth Riemannian metric $m$ on $\Ss^2$ (\ie every
element of $G$ is an isometry of $m$); and then, by the Uniformization
theorem, there is a diffeomorphism $h\colon\Ss^2\carr$ such that the
pull-back metric $h^\star(m)$ is conformally equivalent to the
standard metric $m_0$ on $\Ss^2$. Consequently, $h^{-1}Gh<\SL_2(\C)$,
the group of Möbious transformations of $\Ss^2$. Since every element
of $G$ is topologically conjugates to finite order rotations, this
implies $h^{-1}Gh<\SO(3)$.

In order to construct our $G$-invariant Riemannian metric $m$, we
start considering an arbitrary metric $m'$ on $\Ss^2$. For each
$\e>0$, we define the Riemannian metric
\begin{displaymath}
  m^{\e} :=  \sum_{g \in G} e^{-\e|g|_S}g^* m',
\end{displaymath}
where $S\subset G$ is a finite set of generators of $G$ and
$\abs{\cdot}_S$ denotes the associated word length function on
$G$. Due to the sub-exponential growth of the group $G$ and the
sub-exponential growth of derivatives, we show these metrics are
well-defined and smooth, for every $\e>0$. Then, we show in
Lemma~\ref{metrics} that each element of the generating set
$S\subset G$ is $e^\e$-quasiconfomal with respect to $m^\e$. Then,
invoking rather classical arguments about compactness of the space of
quasiconformal maps (Lemma~\ref{quasi}), we get our $G$-invariant
Riemannian metric $m$ as a limit of metrics $m^\e$, with $\e\to 0$,
which completes the proof of Theorem~\ref{thm:main}.

Finally, in \S\ref{sec:burnside-problem-Td-hyperbolic-mfds} we prove
Theorems~\ref{thm:Burnside-problem-hyperbolic-mfd} and
\ref{thm:Burnside-problem-Td}.

\subsection{Acknowledgments}
\label{sec:acknowledgments}

We want to thank Clark Butler, Amie Wilkinson and Aaron Brown for
useful comments and to Ian Agol for encouragement to the first named
author. A.K. was partially supported by FAPERJ-Brazil and CNPq-Brazil.

\section{Preliminaries}
\label{sec:preliminaries}

In this section we fix some notation we shall use all along the paper
and recall some concepts and results. 

\subsection{Finitely generated groups}
\label{sec:finit-gener-groups}

Let $G$ be a finitely generated group and $S\subset G$ be a finite set
of generators of $G$. We say $S$ is \emph{symmetric} when $s^{-1}\in
S$, for every $s\in S$.

Then given a symmetric set of generators $S$, we define the \emph{word
  length function} $\abs{\cdot}_S\colon G\to\N_0$ by $\abs{id}_S=0$
and
\begin{equation}
  \label{eq:word-length-def}
  \abs{g}_S:=\min\{n\in\N : g=s_{j_1}s_{j_2}\ldots s_{j_n},\ \text{for }
  s_{j_i}\in S\}, 
\end{equation}
for every $g\in G\setminus\{id\}$.

We say that the group $G$ has \emph{sub-exponential growth} when it
holds
\begin{equation}
  \label{eq:sub-exponential-growth}
  \lim_{n\to+\infty} \frac{\log\sharp\left\{ g\in G : \abs{g}_S\leq
      n\right\}}{n}  = 0.
\end{equation}
where $\sharp\{\cdot\}$ denotes the amount of elements of the set.  It
is well known that this concept does not depend on the finite set of
generators.

Observe that, by classical sub-additive arguments, the above limit
\eqref{eq:sub-exponential-growth} always exists.

\subsection{Groups of diffeomorphisms and $C^r$-norms}
\label{sec:diffeos-Cr-norms}

Let $M$ be a closed smooth manifold. The group of
$C^r$-diffeomorphisms of $M$ will be denoted by $\Diff{r}{}(M)$. The
sub-group of $C^r$-diffeomorphisms which are isotopic to the identity
shall be denoted by $\Diff{r}0(M)$.

When $M$ is orientable and $\omega$ is a smooth volume form on $M$, we
write $\Diff{r}\omega(M)$ for the group of $C^r$-diffeomorphisms that
leaves $\omega$ invariant.

When $M$ is endowed with a Riemannian metric $m$, we write $UTM$ for
the \emph{unit tangent bundle,} \ie
\begin{displaymath}
  UTM:=\{v\in TM : \abs{v})=1\}, 
\end{displaymath}
where $\abs{\cdot}$ denotes the norm induced by $m$.

Then, given any $f\in\Diff1{}(M)$ we define
\begin{equation}
  \label{eq:D+f-def}
  \norm{D(f)}^+:=\sup_{v\in UTM} \abs{Df(v)}, 
\end{equation}
and
\begin{equation}
  \label{eq:D-f-def}
  \norm{D(f)}:=\max\left\{\norm{D(f)}^+,\norm{D(f^{-1})}^+\right\}. 
\end{equation}

Now, we will define the $C^r$-norm $\norm{\cdot}_r$ (see
\cite[\S~4]{FisherMargulisAlmostIsom} and the references therein for a
more complete discussion).

Let us start noticing that the Riemannian metric $m$ induces an
isomorphism (by duality) between $T(M)$ and $T^*(M)$, \ie the tangent
and cotangent bundles of $M$. So, $g$ induces an inner product on both
bundles, and consequently, on every tensor bundle over $M$. For
instance, $g$ induces a metric in the space of symmetric $i$-tensors,
which is denoted by $S^i(T^{*}(M))$, for any $i$. So, we shall
consider the space of Riemannian metrics endowed with the topology
induced by the ambient space $S^2(T^{*}(M))$.

The norms $\norm{\cdot}_r$ are defined using the language of
jets. Given a $C^{r}$ bundle $E$ over $M$, let $J_r(E)$ be the vector
bundle of $r$-jets on $E$ (for more details, see
\cite{FisherMargulisAlmostIsom}, \cite[Chapter
1]{EliashbergMishachevHPrinc} and references therein). So, a $C^{r}$
section of $E$ gives a continuous section of $J_r(E)$ (but not the
other way around) and two such sections coincides at some point in $M$
if and only if the derivatives of original sections agree up to order
$r$ at that point.

Observe there is a natural identification
\begin{displaymath}
  J_r(E) \cong \bigoplus_{i=1}^r
  S^i(T^{*}M)\otimes E
\end{displaymath}
as proven, for example, in \cite[\S 4]{FisherMargulisAlmostIsom}.

Let $J^r(M)$ be the bundle of $r$-jets of sections of the trivial
bundle $E := M \times \R$. Then notice that any $C^r$-diffeomorphism
$\phi\colon M\carr$ naturally induces a linear map
$j^r(\phi)(x)\colon J^r(M)_x \to J^r(M)_{\phi^{-1}(x)}$ that sends
each $C^r$-section $s\colon M \to \R$ to $s \circ \phi$. Therefore, we
can define
\begin{equation}
  \label{eq:Cr-norm-def}
  \norm{\phi}_r:= \max_{x \in M}\norm{j^r(\phi)(x)},
\end{equation}
where $\norm{j^r(\phi)(x)}$ is the operator norm defined from the
norms on the vector spaces $J^r(M)_x$ and $J^r(M)_{\phi^{-1}(x)}$.

One can also define norms $\norm{\cdot}_r'$ on the bundles $J^r(M)$
using coordinate charts as follows: consider a finite covering of $M$
by coordinate charts $(U_i,\psi_i)$. For a $C^{r}$ section $s$ of
$E = M \times\R$, one defines
\begin{displaymath}
  \norm{s}_r' := \max_{x, j, i } \norm{D_{\psi_i(x)}^j (s\circ
    \psi_i^{-1})},
\end{displaymath}
where this maximum is taken considering the functions
$s\circ \psi_i^{-1}\colon\R^n \to \R$ in each coordinate chart, and
then calculating the maximum absolute value of a partial derivative of
degree $j$ less or equal to $r$ of such functions over $x \in U_i$.

It is easy to verify that the norms $\norm{\cdot}'_r$ and
$\norm{\cdot}'_r$ are equivalent on $J^r(M)$.

Given a diffeomorphism $\phi\in\Diff{r}{}(M)$, similarly as we define
the norm $\norm{\phi}_r$, we can define a norm $\norm{\phi}'_r$ and
these norms must be equivalent.

We will need the following basic fact relating $\norm{\cdot}_1$ and
$\norm{D(\cdot)}$ we defined in \eqref{eq:D-f-def}:

\begin{proposition}
  \label{tec1}
  Let $M$ be a closed smooth manifold. Then there exists a constant
  $C>0$ such that 
  \begin{displaymath}
    \norm{\phi}_1 \leq C \norm{D(\phi)},
    \quad\forall\phi\in\Diff1{}(M). 
  \end{displaymath}
\end{proposition}

\begin{proof}
  This follows from the fact that the norms $\norm{\cdot}_1'$ and
  $\norm{\cdot}_1$ are equivalent, $\norm{\phi}_1'$ is bounded above
  by the maximum derivative of $\phi$ in coordinate charts, which is
  also bounded by a $C\norm{D(\phi)}$, for some constant $C$ just
  depending on $M$.
\end{proof}

The following notion plays a fundamental role in our work:

\begin{definition}[Sub-exponential growth of derivatives]
  \label{def:sub-expon-growth-deriv}
  Let $M$ be a smooth closed manifold, $G<\Diff\infty{}(M)$ be a
  finitely generated subgroup and $S\subset G$ a finite symmetric set
  of generators. Then we say that $G$ has \emph{sub-exponential growth
    of derivatives} when for every $\e>0$ and any $r\in\N$, there
  exists $N_{\e,r}\in\N$ such that
  \begin{displaymath}
    \norm{w_n}_r\leq e^{n\e}, \quad\forall n\geq N_{\e,r},
  \end{displaymath}
  and any $w_n\in G$ such that $\abs{w_n}_S\leq n$, where
  $\abs{\cdot}_S$ denotes the word length function given by
  \eqref{eq:word-length-def}. 
\end{definition}

\subsection{Riemannian metrics}
\label{sec:riemannian-metrics}

In \S~\ref{sec:sub-expon-growth-case} we shall need to deal with
convergence of Riemannian metrics, so here we recall some basic facts
about the space of metrics.

Recall that if $M$ is a closed smooth manifold, a Riemannian structure
$m$ on $M$ induces a metric on $S^2(T^{*}(M))$. If $s$ denotes a
section of the tensor bundle $S^2(T^{*}(M))$ and
$\phi \in \Diff{r+1}{}(M)$, there is a natural action of $\phi$ in the
bundle of $r$-jets of sections of $S^2(T^{*}(M))$ which is defined by
sending $s$ to the pullback ${\phi}^{*}(s)$. We will need the
following basic

\begin{proposition}
  \label{tec2}
  There exists a constant $C>0$ just depending on $M$ such that
  \begin{displaymath}
    \norm{\phi^{*}(s)}_{r} \leq C\norm{\phi}_{r+1}\norm{s}_r. 
  \end{displaymath}
\end{proposition}

\begin{proof}
  By taking coordinate charts, we can also define norms on
  $J^r(S^2(T^{*}(M)))$ as we did for $J_r(M)$, coinciding with the
  norms for metrics defined in
  \cite{HurtadoContDiscHomoDiffeoGr}. Then, Proposition \ref{tec2}
  follows from Lemma 3.2 in \cite{HurtadoContDiscHomoDiffeoGr}.
\end{proof}

We will also need to make use of the following

\begin{proposition}
  \label{tec3}
  Let $r \geq 1$ and $(s_n)_n$ be a Cauchy sequence in the space of
  $r$-jets of $S^2(T^{*}(M))$. Then there exists a section $s$ of the
  bundle of $r$-jets in $S^2(T^{*}(M))$ such that
  $\norm{s_n - s}_{r} \to 0$, as $n\to+\infty$.
\end{proposition}

\begin{proof}
  This is consequence of the definition of the norms $\norm{\cdot}_r'$
  in terms of coordinate charts , the fact that a Cauchy sequence
  of $C^r$ real functions which are supported in a compact set
  converges to a $C^{r}$ function real function in the $C^{r}$
  topology, which is an easy consequence of Arzela-Ascoli theorem.

  More details can be found in \cite{HurtadoContDiscHomoDiffeoGr}.
\end{proof}

\section{Sub-exponential growth of derivatives}
\label{sec:derivatives}

All along this section, $M$ will denote a closed orientable surface
and $\omega$ a smooth area form on $M$.

Given any $f\in\Diff1{}(M)$, we say $p\in M$ is a \emph{hyperbolic
  fixed point} of $f$ when $f(p)=p$ and the spectrum of $Df_p$ does
not contains complex numbers of modulus equal to $1$. 

\begin{definition}
  \label{def:ellipti-group}
  We say that a group $G<\Diff1{}(M)$ is \emph{elliptic}, if there is
  no element in $G$ having a hyperbolic fixed point.
\end{definition}

For example, any subgroup of $\SO(3)$ is an elliptic subgroup of
$\Diff\infty\omega(\Ss^2)$, where $\omega$ denotes the smooth area
form induced from the Euclidean structure of $\R^3$. There exist other
examples of elliptic sub-groups of $\Diff\infty\omega(\Ss^2)$. For
example, one can construct abelian groups of commuting
pseudo-rotations using the so called Anosov-Katok method (see for
instance \cite{FayadKatokConstElliptic}).

A natural question about elliptic groups of diffeomorphisms of $\Ss^2$
is the following:

\begin{question}
  \label{ques:elliptic-groups-S2}
  Are all elliptic subgroups of $\Diff\infty\omega(\Ss^2)$ either
  solvable or conjugate to a subgroup of $\SO(3)$?
\end{question}

The main result of this section is the following:

\begin{lemma}[Sub-exponential growth of derivatives]
  \label{ideal}
  Let $M$ be a closed orientable surface, $\omega$ an area form on $M$
  and $G$ be a finitely generated elliptic sub-group of
  $\Diff\infty\omega(M)$. Then, $G$ has sub-exponential growth of
  derivatives (see Definition~\ref{def:sub-expon-growth-deriv}).
\end{lemma}

We will begin the proof of Lemma~\ref{ideal} by proving a similar
weaker result that just consider the first derivative. In fact, we
will start proving the following 

\begin{lemma}
  \label{ideal2}
  Let $G$ be a finitely generated elliptic subgroup of
  $\Diff{1+\alpha}\omega(M)$, for some $\alpha>0$, and let $S$ be a
  finite symmetric generating set of $G$. Then, for any $\e>0$, there
  exists $N_{\e} > 0$ such that 
  \begin{displaymath}
    \norm{D(w_n)}\leq e^{\e n}, \quad\forall n\geq N_\e,
  \end{displaymath}
  and every $w_n \in G$ with $\abs{w_n}\leq n$.
\end{lemma}

To prove Lemma~\ref{ideal2}, we first need some definitions. Let
$S = \{s_1,s_2,...,s_r\}$ be the generating symmetric set of $G$ and
consider the space $\Sigma := S^\Z$ consisting of bi-infinite
sequences of elements of $S$.  There is a natural \emph{shift map}
$\sigma\colon\Sigma\carr$ given by
$\sigma : (\ldots, g_{-1}, g_{0}, g_1,\ldots) \mapsto (\ldots,
g_{-1}', g_{0}', g_1',\ldots)$, where $g_i':=g_{i+1}$, for every
$i\in\Z$. Then consider the map $F\colon\Sigma \times M\carr$ given by
\begin{displaymath}
  F(w,x)= \big(\sigma(w), g_0(x)\big), \quad\forall
  w=(\ldots,g_1,g_0,g_1\ldots)\in\Sigma,\ \forall x\in M.
\end{displaymath}
The map $F$ encodes the group action and we have the obvious
commutative diagram:

\begin{displaymath}
  \begin{CD}
    \Sigma \times M @>F>>  \Sigma \times M \\
    @V{\pi}VV          @VV{\pi}V \\
    \Sigma @>>{\sigma}> \Sigma
  \end{CD}
\end{displaymath}

The idea of the proof of Lemma~\ref{ideal2} goes as follows: reasoning
by contradiction, we suppose sub-exponential growth of derivatives
does not hold. Then, we show the existence of an ergodic $F$-invariant
probability measure $\mu$ exhibiting non-zero top Lyapunov exponent.
along the fiber, \ie on $M$. Since we are assuming the fiber $M$ has
dimension $2$ and the action preserves area\footnote{In fact, this is
  the only point where the volume preserving assumption in Theorem
  \ref{thm:main} is crucial.}, then the lower Lyapunov exponent must
be negative. On the other hand, the dynamics on the base, given by the
full-shift, is uniformly hyperbolic. So, $\mu$ is essentially a
hyperbolic measure (see \cite[page 659]{KatokHasselblatt} for
details).

We now recall the following result due to Katok:

\begin{theorem}[Katok~\cite{KatokLyapExp,KatokHasselblatt}]
  \label{Katok}
  Let $N$ be a compact manifold, $F$ be a $C^{1 + \alpha}$
  diffeomorphism of $N$ and $\mu$ be an ergodic hyperbolic measure for
  $F$. Then, there exists a hyperbolic periodic point of
  $F$. Moreover, the periodic point (and its orbit) can be chosen as
  close to $\text{supp}(\mu)$ as one wants.
\end{theorem}

We will prove that Theorem~\ref{Katok} is also true if one considers
the space $N := \Sigma \times M$ and the map $F$ as before (observe
that $\Sigma \times M$ is not a manifold and $F$ is not a
diffeomorphism) and so we will obtain a hyperbolic periodic point
$(w,x)$ for $F$ of say order $k$. The element $w \in \Sigma$ is
determine by the infinite bi-concatenation of a word of length $k$
that defines an element of $G$ having $x$ as a hyperbolic fixed point,
getting a contradiction.

To avoid reproving Katok's theorem for our space $\Sigma \times M$ and
$F$, we will embed the dynamics of $F$ into the dynamics of a
diffeomorphism $F'$ of a 4-manifold $N$. Then, we will construct a
hyperbolic measure $\mu$ for $F'$ and directly apply Katok's theorem
as stated above to $F'$.

Now, we prove Lemma~\ref{ideal2}:

\begin{proof}[Proof of Lemma \ref{ideal2}] First we embed the dynamics
  of the shift $\sigma\colon\Sigma\carr$ into the dynamics of a linear
  $2$-dimensional horseshoe map (see \cite[Chapter 2, \S
  5]{KatokHasselblatt} for details). More formally, there exists a
  $C^\infty$-diffeomorphism $h\colon\Ss^2\carr$ such that $h$ acts on
  an open set $U \subset\Ss^2$ as the linear Smale's horseshoe map and
  where $\Lambda:=\bigcap_{n \in \Z} h^n(U) $ is a hyperbolic Cantor
  set for $h$.  So, there is an embedding $E'\colon\Sigma \to\Ss^2$
  such that $E'(\Sigma) = \Lambda$ and the following commutative
  diagram holds:

  \[
    \begin{CD}
      \Sigma @>E'>>  \Ss^2 \\
      @V{\sigma}VV          @VV{h}V \\
      \Sigma @>>E'> \Ss^2
    \end{CD}
  \]

  The embedding $E'$ can be naturally extended to an embedding
  $E\colon\Sigma \times M \to\Ss^2\times M$. Then we will extend the
  homeomorphism
  $E\circ F\circ E^{-1}\colon E(\Sigma \times M) \to E(\Sigma \times
  M)$ to a smooth diffeomorphism $F'\colon\Ss^2\times M\carr$ which is
  a skew-product over the diffeomorphism $h\colon\Ss^2\carr$. This
  extension exists only if $G \subset \Diff{2}{0}(M)$, \ie every
  element of $G$ is isotopic to the identity. For the sake of
  simplicity, and since we are mainly interested in the case
  $M=\Ss^2$, we will assume $G$ is contained in the identity isotopy
  class.

  Then such an extension is constructed considering a smooth map
  $f\colon\Ss^2\carr\Diff20(M)$ (\ie
  $\Ss^2\times M\ni (z,x)\mapsto f_z(x)\in M$ is $C^2$) such that
  $f_z=g_0$, for every $z\in\Lambda\subset\Ss^2$, and
  $(\ldots,g_{-1},g_1,g_1,\ldots)=E^{-1}(z)\in\Sigma$. In fact, given
  such a map $f$, we can simply define
  \begin{displaymath}
    F'(z,x):=\big(h(z),f_z(x)\big), \quad\forall (z,x)\in\Ss^2\times
    M,
  \end{displaymath}
  and clearly get the following commutative diagram: 
  \begin{displaymath}
    \begin{CD}
      \Sigma \times M @> E >>  \Ss^2 \times M \\
      @V{F}VV          @VV{F'}V \\
      \Sigma \times M @>>E> \Ss^2 \times M
    \end{CD}
  \end{displaymath}
    
  In conclusion, we have embedded the dynamics of $F$ into the
  dynamics of a $C^2$-diffeomorphism $F'$ of the $4$-manifold
  $N := \Ss^2\times M$ that fibers over a diffeomorphism
  $h\colon\Ss^2\carr$. We will now construct the ergodic hyperbolic
  measure $\mu$ for $F'$.

  Let $UTM$ be the unit tangent bundle of $M$. Then we consider the
  map $\partial F'\colon\Ss^2\times UTM\carr$ that fibers over $F'$
  which is given by
  \begin{displaymath}
    \partial F'\big(z,(x,v)\big):=\left(h(z),
      \frac{Df_z(x)v}{\norm{Df_z(x)v}}\right), \quad\forall
    z\in\Ss^2,\ \forall (x,v)\in UTM.
  \end{displaymath}

  Analogously, one can define the map $\partial F\colon\Sigma\times
  UTM\carr$ that fibers over $F$ and is given by
  \begin{displaymath}
    \partial F\big(w,(x,v)\big):=\left(\sigma(w),
      \frac{Dg_0(x)v}{\norm{Dg_0(x)v}}\right), 
  \end{displaymath}
  for every $w=(\ldots,g_{-1},g_0,g_1,\ldots)\in\Sigma$ and every
  $(x,v)\in UTM$; and the map $\partial E\colon\Sigma\times UTM\carr$
  just given by $\partial E:=E'\times id_{UTM}$.

  Notice that by our definitions, it holds
  $\partial E \circ \partial F = \partial F'\circ\partial E$.

  Then, let us suppose there is no sub-exponential growth of
  derivatives. So, there exists $\e>0$, a sequence of words $w_n$ of
  elements of $S$ of length smaller or equal than $n$, and vectors
  $(x_n,v_n) \in UTM$ such that $\norm{D_{x_n}w_n(v_n)}\geq e^{\e n}$,
  for each $n\geq 1$.  Let us define ${\hat{w}_n} \in \Sigma$ to be
  the bi-infinite periodic word ${\hat{w}_n} =\ldots w_nw_nw_n\ldots$
  and let $\delta_{({\hat{w}_n}, (x_n, v_n))}$ be the Dirac measure on
  $\Sigma\times UTM$ supported on the point
  $({\hat{w}_n}, (x_n, v_n))$. We consider the sequence of measures
  \begin{displaymath}
    \nu_n := \frac{1}{n} \sum_{i=1}^n
    \big((\partial F)^i\big)_{*}(\delta_{({\hat{w}_n}, (x_n, v_n))}),
    \quad\forall n\geq 1.
  \end{displaymath}
  
  Then, consider the function $\psi\colon\Ss^2\times UTM \to \R$ given
  by
  \begin{displaymath}
    \psi\big(z,(x,v)\big) = \log
    \bigg(\frac{\norm{Df_z(v)}_{f_z(x)}}{\norm{v}_{z}}
    \bigg), \quad\forall z\in\Ss^2,\ \forall (x,v)\in UTM,
  \end{displaymath}
  and observe that
  \begin{equation}
    \label{birk}
    \begin{split}
      \int_{\Ss^2\times UTM} \psi\dd(\partial E_*\nu_n) &=
      \int_{\Sigma\times UTM} (\psi\circ\partial
      E)\dd\nu_n \\
      &= \frac{1}{n}\sum_{i=1}^n \psi\circ \partial E\Big(\partial
      F^i\big(\hat w_n, (x_n,v_n)\big)\Big) \\
      &\geq
      \frac{1}{n}\log\left(\frac{\norm{D_{x_n}w_n(v_n)}_{w_n(x_n)}}
        {\norm{v_n}_{x_n}}\right)\geq\e,
    \end{split}
  \end{equation}
  for every $n\geq 1$.

  So, by Banach-Alaoglu theorem, there exists a sub-sequence
  $(\partial E_*\nu_{n_j})_{n_j}$ that converges in the weak-star
  topology to a measure $\nu'$, which is clearly
  $\partial F'$-invariant, and by \eqref{birk}, it holds
  \begin{displaymath}
    \int_{\Ss^2\times UTM}\psi\dd\nu'\geq\e. 
  \end{displaymath}

  Then, if we consider the ergodic decomposition of $\nu'$, there
  exists an ergodic $\partial F'$-invariant probability measure $\nu$
  such that
  \begin{equation}
    \label{eq:Lyap-exp-cond-inv-meas}
    \int_{\Ss^2\times UTM} \psi\dd\nu \geq \e.
  \end{equation}

  Now, let $\mu$ be the push-forward measure of $\nu$ by the
  projection on the $\Ss^2\times M$ factor, \ie
  $\mu:={\mathrm{pr}}_*\nu$, where
  $\mathrm{pr}\colon\Ss^2\times UTM\to\Ss^2\times M$ denotes the
  natural projection. Observe $\mu$ is an ergodic $F'$-invariant
  measure. We claim $\mu$ is a hyperbolic measure for $F'$, \ie all
  its Lyapunov exponents are different from zero.

  In order to prove that, first observe that measure $\mu$ is
  supported on the subset $\Lambda\times M$, where $\Lambda$ is the
  horseshoe of diffeomorphism $h$, and consequently, it is a uniform
  hyperbolic set. So, to show that $\mu$ is a hyperbolic measure, it
  is enough to show that its both Lyapunov exponents along the
  vertical fibers (\ie on the $M$ factor) are different from zero.

  Now, combining \eqref{eq:Lyap-exp-cond-inv-meas} and a result due to
  Ledrappier~\cite[Proposition 5.1]{LedrappierQuelProp} one conclude
  the the top Lyapunov exponent of $\mu$ along vertical fibers is
  positive. Since the diffeomorphism $f_z\in\Diff2{0}(M)$ leaves
  invariant the area form $\omega$, for every $z\in\Lambda$, and $\mu$
  is supported on $\Lambda\times M$, this implies the bottom Lyapunov
  exponent of $\mu$ along vertical fibers is negative. Thus, $\mu$ is
  an ergodic hyperbolic measure.

  So, we can apply Theorem~\ref{Katok} to conclude that, for every
  open neighborhood $V$ of $\Lambda$, there is a hyperbolic periodic
  point $p'\in V\times M$ for $F'$. Since $\Lambda$ is a locally
  maximal invariant set for $h$, \ie $\Lambda=\bigcap_{n\in\Z} h^n(V)$
  for every sufficiently small neighborhood $V$ of $\Lambda$, this
  implies $p'\in\Lambda\times M$. Hence we can consider the point
  $(w^p,x^p):=E^{-1}(p')\in\Sigma\times M$, which is a periodic point for
  $F$. Since $\partial F'\circ\partial E=\partial E\circ\partial
  F$, we conclude $(w^p,x^p)$ is hyperbolic along the vertical fibers,
  contradicting the fact that there is a natural number $k$ such that
  $F^k(w^p,x)=(w^p,x)$, for every $x\in M$.
\end{proof}

In order to finish the proof of Lemma \ref{ideal}, first we need to
recall the following result which is an easy consequence of the chain
rule in dimension one and, according to Fisher and
Margulis~\cite{FisherMargulisAlmostIsom}, it is a well known estimate
used in KAM-theory:

\begin{lemma}[Lemma 6.4 in \cite{FisherMargulisAlmostIsom}]
  \label{marfish}Let $M$ be a compact smooth manifold, let
  $\phi_1, \phi_2,...\phi_n \in \Diff{r}{}(M)$,
  $N_1 := \max_{1\leq i \leq n} \|\phi_i\|_1$ and
  $N_r := \max_{1\leq i \leq n} \|\phi_i\|_r$. Then, there exists a
  polynomial $Q$ (just depending on the dimension of $M$) such that
  \begin{displaymath}
    \norm{\phi_1 \circ \phi_2.....\circ \phi_n}_r \leq
    N_1^{rn}Q(nN_r).  
  \end{displaymath}
\end{lemma}

We can now give a proof of Lemma \ref{ideal}.

\begin{proof}[Proof of Lemma~\ref{ideal}]
  From Proposition \ref{tec1} and Lemma \ref{ideal2}, one observes
  that Lemma \ref{ideal} holds for $r = 1$. Then, given any $r > 1$
  and $\e>0$, we apply Lemma \ref{ideal} to guarantee the existence of
  a natural number $N_{\frac{\e}{3r},1}$ so that
  \begin{displaymath}
    \norm{w_n}_1 \leq e^{\frac{\e n}{3r}}, \quad\forall n \geq
    N_{\frac{\e}{3r},1},
  \end{displaymath}
  and any $C r$-diffeomorphism $w_n\colon M\carr$ that can be written
  as a composition of at most $n$ elements of $S$.
  
  Then, by Lemma~\ref{ideal} easily follows from Lemma \ref{marfish}
  and a proper chosen of constant $N_{\e,r}$.
\end{proof}

\section{The sub-exponential growth case}
\label{sec:sub-expon-growth-case}

In this section we prove Theorem~\ref{thm:main} for periodic groups of
sub-exponential growth. So, we will assume $G<\Diff{\infty}{}(\Ss^2)$
is a finitely generated periodic sub-group with bounded exponent and
sub-exponential growth, \ie if $S\subset G$ is a finite generating set
of $G$, then condition \eqref{eq:sub-exponential-growth} holds.

\begin{theorem}
  \label{subexcase}
  Let $\omega$ be an area form on $\Ss^2$ and
  $G<\Diff{\infty}\omega(\Ss^2)$ be a finitely generated periodic
  sub-group with bounded exponent and sub-exponential growth (\ie if
  $S\subset G$ is a finite generating set of $G$, then condition
  \eqref{eq:sub-exponential-growth} holds).

  Then, $G$ is finite. 
\end{theorem}

The idea of the proof of Theorem~\ref{subexcase} goes as follows:
combing the sub-exponential growth of $G$ and the sub-exponential
growth of derivatives we proved in Lemma~\ref{ideal} we find, for each
$\e>0$, a $C^{\infty}$ \emph{``à la Pesin''} Riemannian metric
$m^{\e}$ on $\Ss^2$ such that each element of $S$ is
$e^{\e}$-bi-Lipschitz for $m^{\e}$ (Lemma \ref{metrics}). By the
Uniformization theorem of surfaces, each metric $m^{\e}$ is
conformally equivalent to the standard metric $m_0$ in $\Ss^2$, which
implies there are conjugates $G_{\e}= g_{\e}^{-1}Gg_{\e}$ of the group
$G$ such that each element of the generating set
$S_{\e}:= g_{\e}^{-1}Sg_{\e}$ is an $e^{\e}$-quasi-conformal
homeomorphism of $\Ss^2$. Then, as $\e \to 0$, we will attempt to show
that $G_{\e}$ converges to a Burnside group $G_0$ which acts
conformally on $\Ss^2$, \ie the group $G_0$ is a periodic subgroup of
$\PSL_2(\mathbb{C})$ which, by Schur theorem, must be finite. We will
show this implies $G$ is finite as well.

We start constructing the family of Riemannian metrics $(m^\e)_{\e>0}$
which is mainly motivated by Pesin's work on non-uniform hyperbolicity
(see for instance \cite{KatokHasselblatt} for details):

\begin{lemma}
  \label{metrics}
  Let $M$ be a closed smooth manifold and $G<\Diff\infty{}(M)$ be a
  finitely generated group with sub-exponential growth and
  sub-exponential growth of derivatives (see
  Definition~\ref{def:sub-expon-growth-deriv}), and $S$ be a finite
  set of generators of $G$.  Then for any $\e>0$ there exists a
  $C^{\infty}$-metric $m^{\e}$ such that
  \begin{equation}
    \label{eq:m-e-Pesin-like}
    e^{-\e}\abs{v}_{x}^{\e} \leq \abs{D_xs(v)}_{s(x)}^{\e} \leq
    e^{\e}\abs{v}_{x}^{\e}, 
  \end{equation}
  for every $s \in S$ and every $(x,v) \in TM$, and where
  $\abs{v}^\e:=\sqrt{m^\e_x(v,v)}$. 
\end{lemma}

\begin{proof}

  Let us start considering an arbitrary $C^\infty$ Riemannian metric
  $m$ on $M$ and $\e$ be any positive number.  For each integer $n>0$,
  let us write
  \begin{displaymath}
    m^{\e}_n =  \sum_{g \in B_n} e^{-\e\abs{g}_S}g^{*}m,
  \end{displaymath}
  where $B_n:=\{g\in G : \abs{g}_S\leq n\}$ and $g^*m$ denotes the
  pull-back metric of $m$ by $g$, \ie
  $g^*m(v,w):=m\big(Dg(v),Dg(w)\big)$.  

  The metric $m^\e$ will be constructed as the limit of the sequence
  $\big(m^{\e}_n\big)_n$, as $n\to\infty$.

  Fixing an integer $r>0$, by Lemma \ref{ideal} there is a natural
  number $N_{\e/2,{r+1}} >0 $ so that for any $g \in G$ satisfying
  $\abs{g}_S\geq N_{\e/3,{r+1}}$, we have
  \begin{displaymath}
    \norm{g}_{r+1} \leq e^{\frac{\e\abs{g}_S}{3}}.
  \end{displaymath}

  On the other hand, since we are assuming $G$ is a group of
  sub-exponential growth, there exists a constant $K_{\e}> 0$ so that
  \begin{displaymath}
    \sharp\left\{g\in G : \abs{g}_S\leq n\right\}  \leq e^{\frac{\e
        n}{3}}, \quad\forall n\geq K_\e. 
  \end{displaymath}

  Therefore for any $n \geq \max \{{N_{\e/3,r}, M_{\e}} \}$, by
  Proposition \ref{tec3} we have 
  \begin{align*}
    \norm{m^\e_{n} - m^\e_{n-1}}_r &=  \norm{\sum_{g \in B_n \setminus
                                     B_{n-1}} e^{-\e n}g^{*}m}_r \\ 
                                   &\leq \sharp(B_n \setminus B_{n-1})
                                     e^{-\e n}  \max_{g \in B_n
                                     \setminus B_{n-1}}{ \|g^{*}m\|_r}
    \\ 
                                   &\leq \sharp(B_n \setminus B_{n-1})
                                     e^{-\e n}Ce^{\frac{\e
                                     n}{3}}\norm{m}_r\\ 
                                   &\leq C\norm{m}_re^{\frac{-\e
                                     n}{3}}. 
  \end{align*}

  Therefore the sequence $\{m^\e_n\}_{n\geq 1}$, when their elements
  are considered as $r$-jets of $S^2(T^*(M))$, is a Cauchy
  sequence. Then, by Proposition~\ref{tec3}, $\{m^{\e}_n\}_{n\geq 1}$
  converges to a $C^{r}$ tensor $m^\e$ in $S^2(T^*(M))$.

  The tensor $m^{\e}$ is easily seen to be a non-degenerate metric and
  as previous estimates are true for any $r>0$, the metric $m^\e$ is
  in fact $C^{\infty}$, too.
 
  To prove estimates \eqref{eq:m-e-Pesin-like}, observe that given any
  vector $v \in T_x(M)$, we have
  \begin{align*}
    (s^{*}m^{\e}_n)(v,v) &=  \sum_{g \in B_n}
                           e^{-\e|g|_S}(gs)^{*}m(v,v)\\ 
                         &\leq e^{\e}\bigg( \sum_{g \in B_{n+1}}
                           e^{-\e|g|_S}g^{*}m(v,v) \bigg)\\  
                         &= e^{\e} m^{\e}_{n+1}(v,v).
  \end{align*}
  Taking limits in the previous inequality, we obtain
  $\abs{D_xs(v)}_{s(x)}^{\e} \leq e^{\e}\abs{v}_{x}^{\e}$. The other
  inequality follows in a completely analogous way.
\end{proof}

Now, we invoke the Uniformization theorem for the $2$-sphere $\Ss^2$
for finding nice conjugates of our group $G$.  We recall the statement
of the Uniformization theorem in the following form (see \cite[Chapter
10]{DonaldsonRiemannSurf} for more details):

\begin{theorem}[Uniformization of $\Ss^2$]
  \label{thm:uniformization}
  Any $C^{\infty}$ Riemannian metric $m$ on $\Ss^2$ is conformally
  equivalent to the standard metric $m_0$ in $\Ss^2$. That is, there
  exists a $C^{\infty}$ diffeomorphism $g\colon\Ss^2\carr$ and a
  $C^{\infty}$ function $h\colon\Ss^2\to\R$ such that
  \begin{displaymath}
    g^{*}m = e^{h}m_0.
  \end{displaymath}
\end{theorem}

So, invoking Lemma~\ref{metrics} we can construct a family of
Riemannian metrics $(m^\e)_{\e>0}$ that satisfy estimate
\eqref{eq:m-e-Pesin-like} for each element of the generating set
$S$. Then, applying Theorem~\ref{thm:uniformization} to this family of
metrics we can construct a family of $C^\infty$-diffeomorphisms
$(g_{\e})_{\e>0}$ and $C^\infty$ real functions $h_{\e}$ such that
\begin{equation}
  \label{conju}
  g_{\e}^{*}m^{\e} = e^{h_{\e}}m_0, \quad\forall\e>0.
\end{equation}

Now, for each generator $s \in S$ we define
\begin{equation}
  \label{eq:s-e-definition}
  s_{\e} := g^{-1}_{\e}\circ s\circ g_{\e}.
\end{equation}

We will prove that as $\e \to 0$, the diffeomorphisms $s_{\e}$ get
closer to being conformal with respect to the standard metric $m_0$ on
$\Ss^2$. To do that, we first recall the purely metric definition of
quasiconformality (see for instance
\cite{HeinonenKoskelaQuasiconfMapsMet} for more details):

\begin{definition}
  \label{def:quasiconformality}
  Let $(X,d_X), (Y,d_Y)$ be two metric spaces and $f\colon X\to Y$ be
  homeomorphism. For each $r>0$ and $x\in X$ we define
  \begin{displaymath}
    H_f(x,r) := \frac{\sup\{ d_Y\big(f(x),f(y)\big) :  y\in X,\
      d_X(x,y) < r\}}{\inf \{d_Y\big(f(x),f(y)\big) : y\in X,\
      d_X(x,y) > r\}}.
  \end{displaymath}
  Then we say that $f$ is \emph{$K$-quasiconformal}, for some
  $K\geq 1$, whenever
  \begin{displaymath}
    \limsup_{r \to 0} H_f(x,r) \leq K, \quad\forall x \in X.
  \end{displaymath}
\end{definition}

Then we have the following

\begin{lemma}
  \label{lem:s-e-quasiconf}
  For every number $\e >0$ and every $s \in S$, the diffeomorphism
  $s_{\e}\colon\Ss^2\carr$ given by \eqref{eq:s-e-definition} is
  $e^{2\e}$-quasiconformal with respect to the standard Riemannian
  metric $m_0$ on $\Ss^2$. 
\end{lemma}

\begin{proof}
  Let $s$ be an arbitrary element of $S$ and $\e$ a fixed positive
  number.  Given any point $x \in\Ss^2$ and any vector
  $v\in T_x(\Ss^2)$, we write $\abs{v}_x^0$ and $\abs{v}_x^{\e}$ for
  the norms of $v$ in the standard metric $m_0$ and the metric
  $m^{\e}$, respectively. Observe that to prove that $s_{\e}$ is
  $e^{2\e}$- quasiconformal it is enough to show that 
  \begin{equation}
    \label{eq0}
    e^{-2\e}\leq \frac{\abs{D_xs_{\e}(v)}^0_{s_{\e}(x)}}
    {\abs{D_xs_{\e}(w)}^0_{s_{\e}(x)}}\leq e^{2\e}, \quad\forall
    x\in\Ss^2,\  \forall v,w\in T_x(\Ss^2),
  \end{equation}
  such that $\abs{v}_x=\abs{w}_x=1$. 

  Let $y$ be a an arbitrary point in $\Ss^2$ and $v', w'$ be two
  vectors in $T_y(\Ss^2)$. By the inequalities
  \eqref{eq:m-e-Pesin-like} we know that
  \begin{displaymath}
    e^{-\e}\abs{v'}^{\e}_y \leq \abs{D_ys(v')}^{\e}_{s(y)} \leq
    e^{\e}\abs{v'}^{\e}_y.
  \end{displaymath}
  So it follows

  \begin{equation}
    \label{eq1}
    e^{-2\e}\frac{\abs{v'}_y^{\e}}{\abs{w'}_y^{\e}} \leq
    \frac{\abs{D_ys(v')}^\e_{s(y)}}{\abs{D_ys(w')}^\e_{s(y)}} \leq 
    e^{2\e}\frac{\abs{v'}_y^{\e}}{\abs{w'}_y^{\e}}. 
  \end{equation} 

  We also have that as $g_{\e}^{*}m^{\e} = e^{h_{\e}}m_0$. Then,
  $\abs{D_yg_\e(v')}^{\e}_{g_{\e}(y)}= e^{2h_{\e}(y)}\abs{v'}^{0}_y$
  and therefore, if we take $y = g_{\e}(x)$, $v' = D_xg_{\e}(v)$ and
  $w' = D_xg_{\e}(w)$, in inequality \eqref{eq1} we obtain

  \begin{equation}
    \label{eq2}
    \begin{split}
      e^{-2\e} &= e^{-2\e} \frac{e^{2h_{\e}g_{\e}(x)}\abs{v}^{0}_x}
      {e^{2h_{\e}g_{\e}(x)}\abs{w}^{0}_x} \leq
      \frac{\abs{D_{{g_\e}(x)}s(D_xg_{\e}v)}^\e_{sg_{\e}(x)}}
      {\abs{D_{g_\e(x)}s(D_xg_{\e}w)}^\e_{sg_{\e}(x)}} \\
      & \leq e^{2\e}\frac{e^{2h_{\e}g_{\e}(x)}\abs{v}^{0}_x}
      {e^{2h_{\e}g_{\e}(x)}\abs{w}^{0}_x} = e^{2{\e}}.
    \end{split} 
  \end{equation}

  On the other hand, we have that:
  \begin{align*}
    \frac{\abs{D_{{g_\e}(x)}s(D_xg_{\e}v)}^\e_{sg_{\e}(x)}}
    {\abs{D_{g_\e(x)}s(D_xg_{\e}w)}^\e_{sg_{\e}(x)}}
    &= \frac{\abs{D_{s_\e(x)}g_{\e}(D_{x}s_{\e}v)}^\e_{g_{\e}s_{\e}(x)}} 
      {\abs{D_{s_\e(x)}g_{\e}(D_{x}s_{\e}v)}^\e_{g_{\e}s_{\e}(x)}}\\ 
    &= \frac{ e^{2h_{\e}s_{\e}(x)}\abs{D_xs_{\e}v}^{0}_{s_{\e}(x)}}
      {e^{2h_{\e}{s_{\e}(x)}}\abs{D_xs_{\e}w}^{0}_{s_{\e}(x)}}\\  
    &=\frac{\abs{D_xs_{\e}v}^{0}_{s_{\e}(x)}}
      {\abs{D_xs_{\e}w}^{0}_{s_{\e}(x)}}\\  
  \end{align*}

  Therefore, putting together previous inequality and estimate
  \eqref{eq2}, we obtain \eqref{eq0}.
\end{proof}

We will use the following known fact about quasiconformal maps on
surfaces (see for instance \cite[Theorem 1.3.13]{FletcherThesis}):

\begin{lemma}
  \label{quasi}
  Let $x_1,x_2,x_3$ be three different points on $\Ss^{2}$, and
  $(f_{\e}\colon\Ss^2\carr)_{\e>0}$ a family of homeomorphisms
  satisfying:
  \begin{enumerate}
  \item $f_\e$ is $K_\e$-quasiconformal;
  \item $K_{\e} \to 1$, as $\e\to 0$;
  \item for each $\e>0$, $f_{\e}$ fixes $x_1,x_2$ and $x_3$.
  \end{enumerate}

  Then $f_{\e} \to \text{Id}$, as $\e\to 0$, in the $C^{0}$-topology.
\end{lemma}

\begin{remark}
  \label{rem:quasi}
  The same result holds true if the three points $x_i$ are allowed to
  depend on $\e$ under the additional assumption that they remain at a
  bounded distance away from each other.
\end{remark}

\begin{proof}
  This follows from the compactness of $K$-quasiconformal maps of
  $\Ss^2$ which states that the family of $K$-quasiconformal maps of
  $\Ss^2$ fixing three points $x_1,x_2,x_3$ is compact in the
  $C^{0}$-topology (see \cite[Theorem 3.1.13]{FletcherThesis} for
  further details). Therefore, any sub-sequence $\{f_{\e_n}\}$ of maps
  $f_{\e}$ has a convergent sub-sequence converging to a map $f'$ which
  must be $K_{\e}$-quasiconformal for every $K_{\e}$, and therefore,
  conformal. So, $f'\in\PSL_2(\C)$. Besides, such map $f'$ must fix
  $x_1,x_2$ and $x_3$ and so, $f'$ must the identity.
\end{proof}

We will need the following consequence of Lemma \ref{quasi}:
\begin{proposition}
  \label{cantat1}
  Let $s\colon\Ss^2\carr$ be a periodic homeomorphism different form
  the identity, and let $\{s_\e\colon\Ss^2\carr\}_{\e>0}$ be a family of
  homeomorphisms such that the following conditions holds:
  \begin{enumerate}
  \item for each $\e>0$, $s_\e$ is topologically conjugate to $s$;
  \item there is a family of positive real numbers $\{K_\e\}_{\e>0}$
    such that $s_{\e}$ is $K_\e$-quasiconformal, for every $\e>0$ and
    $K_\e \to 1$, as $\e \to 0$;
  \item there exists $\delta > 0$ such that, if $p_\e, q_{\e}$ are the
   two fixed points of $s_\e$, then it holds $d(p_\e,q_\e)>\delta$,
   for every $\e>0$. 
  \end{enumerate}

  Then the family $\{s_{\e}\}_{\e>0}$ is pre-compact with the uniform
  $C^0$-topology and if $s'\colon\Ss^2\carr$ is a homeomorphism such
  that there is a sequence $s_{\e_n}\to s$, with $\e_n\to 0$, then
  $s'\in\PSL_2(\C)$.
\end{proposition}

\begin{proof}
  Given any $\e>0$, the diffeomorphism $s_{\e}$ is conjugate to a
  finite order rotation. So $s_{\e}$ has exactly two fixed points,
  which are denoted by $p_\e$ and $q_\e$. Therefore there exists a
  conformal map $A_\e\in\PSL_2(\C)$ sending $p_{\e}$ to
  $Z:=(0,0,1)\in\Ss^2\subset\R^3$ and $q_{\e}$ to
  $-Z = (0,0,-1)\in\Ss^2\subset\R^3$.  Moreover, we can suppose that
  the family $\{\norm{A_{\e}}\}_{\e>0}$ is bounded and so,
  pre-compact.

  Let us define $S_{\e} := A_{\e}s_{\e}A_{\e}^{-1}$. If we consider
  the great circle $C$ in $\Ss^2$ determined by the $xy$-plane in
  $\R^3$, there is a point $x_{\e} \in C$ such that
  $S_{\e}(x_{\e}) \in C$. Then, there is a rotation
  $R_{\theta_{\e}}\in\SO(3)$ fixing $Z, -Z$ and sending $x_{\e}$ to
  $S_{\e}(x_{\e})$. Therefore, the map $S_{\e} R^{-1}_{\theta_\e}$ is
  still $K_{\e}$-quasiconformal and fixes the three points $Z, -Z $ and
  $S_{\e}(x_{\e})$.

  By Remark~\ref{rem:quasi}, we have $S_{\e}R^{-1}_{\theta_\e} \to Id$
  in the $C^0$-topology, as $\e\to0$. By compactness of $\SO(3)$,
  there is a sub-sequence $\{\e_n\}_n$, with $\e_n\to 0$, as
  $n\to\infty$, and an angle $\theta$ such that
  $R_{\theta_{\e_n}} \to R_{\theta}$, as $n\to\infty$. So,
  $S_{\e_n} \to R_{\theta}$.

  To finish the proof, we must show that $R_{\theta}$ has order $k$.
  As $S_{\e_n} \to R_{\theta}$, as $n\to\infty$, we have that
  $S_{\e_n}^k \to R_{\theta}^k$ and so, $R_{\theta}^k = Id$. This
  implies that $R_{\theta}$ has order $k'$ dividing $k$.

  If $k' < k$, then $S_{\e_n}^{k'} \to Id$, as $n\to\infty$, where
  $S_{\e_n}^{k'}$ is conjugate to a rotation of order
  $\alpha: = \frac{k}{k'}$. But this is impossible as any conjugate of
  a rotation of order $\alpha$ must move some point in $\Ss^2$ at
  least distance $\alpha$ in the standard metric on $\Ss^2$.
\end{proof}

We will now begin the proof of Theorem \ref{subexcase}. We start by
considering the case where $G$ is $2$-generated:

\begin{proposition}
  \label{2gen}
  Let $S$ be a generating set of a group $G$ as in Theorem
  \ref{subexcase}. For any pair of elements $s,t \in S$, the subgroup
  $G' = \langle s,t \rangle$ of $G$ is finite.
\end{proposition}

\begin{proof}
  First let us consider the case where $s,t$ have a common fixed point
  $p$ in $\Ss^2$. The derivative map at $p$ gives a homomorphism
  $\Phi\colon G' \to \GL_2(\R)$. As the image of $\Phi$ is a finitely
  generated periodic linear group, by Schur's theorem such
  homomorphism must have finite image. On the other hand, $\ker(\Phi)$
  is trivial because all the elements of $G'$ are smoothly conjugate
  to rotations. Then, $G'$ must be finite (and in fact, $G'$ must be
  cyclic).

  We will then assume that $s,t$ have no common fixed point.

  Conjugating by the diffeomorphisms $g_{\e}$ given by \eqref{conju},
  we obtain diffeomorphisms $s_{\e}$ and $t_{\e}$ which are
  $e^{2\e}$-quasiconformal. Furthermore, by conjugation with a Möbius
  map we can suppose the fixed points of $s_{\e}$ are $Z:= (0,0,1)$
  and $-Z = (0,0,-1)$ of $\Ss^2 \subset \R^3$ and one fixed points of
  $t_{\e}$ is the point $X = (1,0,0)$.

  Therefore by Proposition \ref{cantat1}, there exists a sequence
  $\{\e_n\}_n$ of positive numbers, with $\e_n\to 0$ as $n\to\infty$,
  so that the sequence $\{s_{\e_n}\}_n$ converges to a non-trivial
  rotation $R_{\theta}\in\SO(3)$, fixing the points $Z$ and $-Z$.
  Recall that each diffeomorphism $t_{\e}$ fixes the point $X$. Let us
  write $X_{\e_n}$ for the other fixed point of $t_{\e_n}$. Then we
  have:

  \begin{lemma}
    \label{confor}
    Either $G'$ is finite, or it holds
    $ \lim_{\e_n \to 0}X_{\e_n} = X$.
  \end{lemma}

  \begin{proof}[Proof of Lemma~\ref{confor}]
    Let us assume the sequence $\{X_{\e_n}\}_n$ does not converge to
    $X$.  So we can invoke Proposition~\ref{cantat1} to guarantee the
    existence of a sub-sequence $\{t_{\e_{n_j}}\}_{n_j}$ converging in
    the $C^{0}$-topology to a Möbius transformation
    $T\in\PSL_2(\C)$. In such a case, the sequences of diffeomorphisms
    $\{s_{\e_{n}}\}_{n}$ and $\{t_{\e_{n_j}}\}_{n_j}$ converge to
    $R_{\theta}$ and $T$, respectively. So, the group
    $\langle R_{\theta}, T \rangle$ is a periodic sub-group of
    $\PSL_2(\mathbb{C})$, and by Schur's theorem, it is finite. Let us
    show this implies $G'$ is finite as well.

    To do that, let $\bb{F}_2$ denote the free group on two elements,
    and let $\{a,b\}\subset\b{F}_2$ be a generating set of
    $\bb{F}_ 2$. Let us write $h\colon\bb{F}_ 2\to G'$ and
    $h_0\colon\bb{F}_ 2\to\PSL_2(\C)$ for the two unique group
    homomorphisms such that $h(a)=s$, $h(b)=t$, and $h_0(a)=R_\theta$,
    $h_0(b)=T$. If $G'$ were infinite, then it would exist an element
    $w\in\bb{F}_ 2$ such that $h(w)\neq id$ and $h_0(w)=id$. However,
    for each $\e_n$ we can consider the only group homomorphism
    $h_{\e_n}\colon\bb{F}_2\to\Diff\infty{}(\Ss^2)$ such that
    $h_{\e_n}(a)=s_{\e_n}$ and $h_{\e_n}(b)=t_{\e_n}$.

    Since each $s_{\e_n}$ and $t_{\e_n}$ is conjugate to $s$ and $t$,
    respectively, we get $h_{\e_n}(w)\neq id$, for every $n$. But on
    the other hand, $h_{\e_n}(w)\to h_0(w)=id$, as $n\to\infty$, which
    contradicts Proposition \ref{cantat1}. So, $G'$ is finite.
  \end{proof}

  We will now deal with the case when $X_{\e_n} \to X$, as
  $n\to\infty$. In order to simplify the notation, we will denote
  $\e_n$ simply by $\e$ and any statement about $\e \to 0$ should be
  understood to be true up to passing to the sequence $\{\e_n\}_n$. We
  will show the following

  \begin{proposition}
    \label{hard}
    If $X_{\e} \to X$, the group $G'$ contains an element of infinite
    order, contradicting the fact that $G'= \langle s,t \rangle$ is
    periodic.
  \end{proposition}

  \begin{proof}[Proof of Propostion~\ref{hard}]
    Since the point $X_\e$ is different of $X$ for every $\e>0$, for
    each $\e$ we can consider the great circle $C_{\e}$ in $\Ss^2$
    passing through the points $X_{\e}$ and $X$. Let
    $M_{\e} \in C_{\e}$ be the midpoint between $X_{\e}$ and $X$, on
    the shortest geodesic segment determined by these points (see
    Figure~\ref{pain2}). Then there is a Möbius transformation
    $A_{\e}$ (a loxodromic element in $\PSL_2(\C)$) such that
    $A_\e(M_\e)=M_\e$, $A_\e(-M_\e)=-M_\e$ and
    $A_\e(X_\e)=-A_\e(X)$. Let us define $Y_\e:=A_\e(X_\e)$. 

    By compactness of $\Ss^2$ and Proposition \ref{cantat1}, there
    exist a point $Y\in\Ss^2$, a rotation $S_\alpha\in\SO(3)$ and a
    sub-sequence of $\{\e_n\}$, such that $Y_{\e} \to Y \in S^{2}$ and
    $A_{\e}t_{\e}A_{\e}^{-1}\to S_\alpha$ is the $C^0$-topology, as
    $\e\to 0$. Notice the rotation $S_\alpha$ has a extricly positive
    finite order and fixes the points $Y$ and $-Y$.

    \begin{figure}[ht]
      \begin{center}
        \includegraphics[scale=0.4]{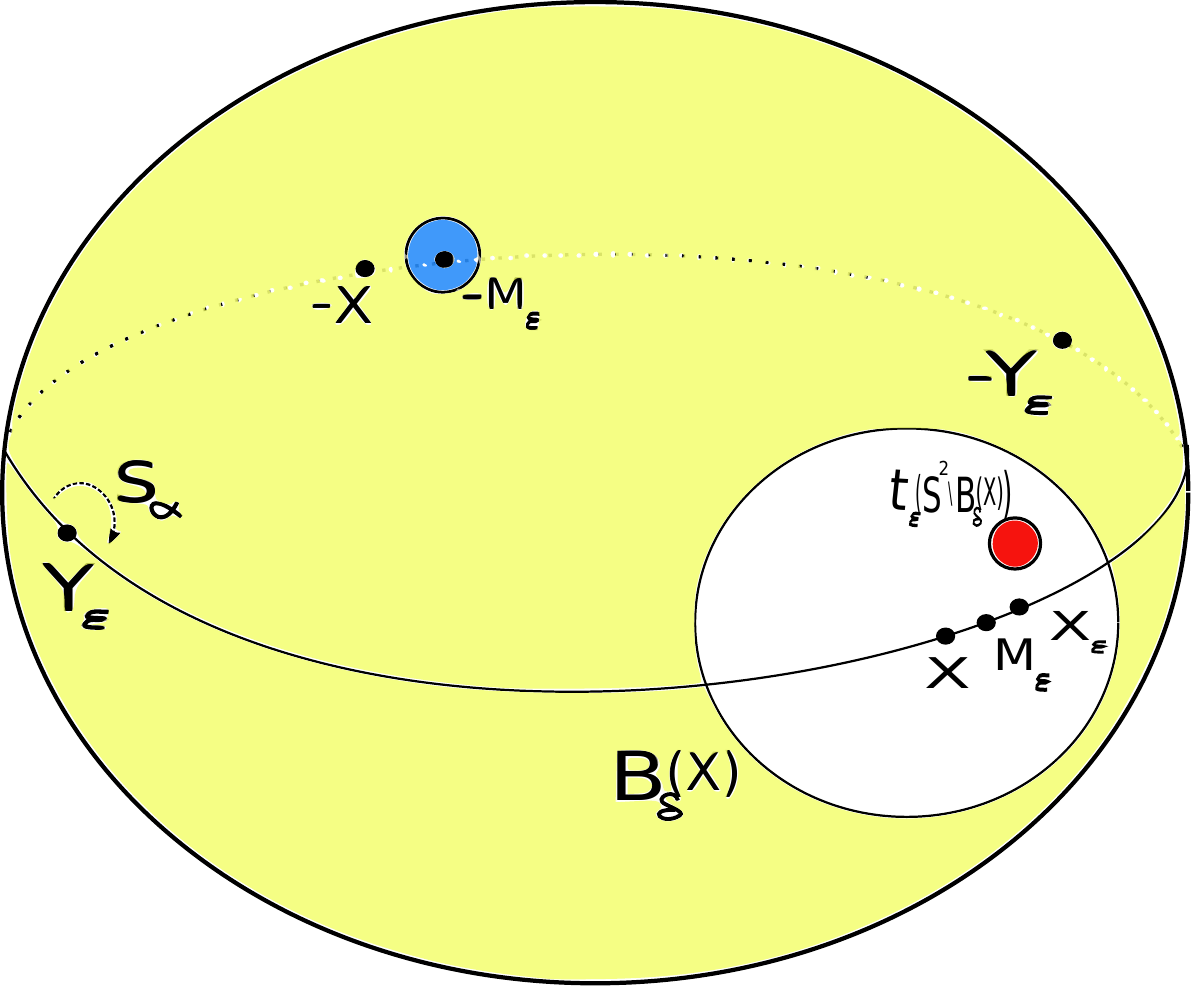}
      \end{center}
      \caption{ The white ball corresponds to $B_\delta(X)$, the blue
        ball to $A_{\e}( \Ss^{2} \setminus B_\delta(X))$ and the red
        ball to $t_{\e}( \Ss^{2} \setminus B_\delta(X))$}
      \label{pain2}
    \end{figure}

    Our purpose now consists constructing an open disk $D\subset\Ss^2$
    such that $s_\e t_\e(D)\subset D$, for $\e>0$ sufficiently small
    enough, and where the inclusion is strict.  This shows the element
    $s_\e t_\e$ cannot have finite order, and this immediately implies
    $st$ has infinite order, as well.
    
    In order to prove that, let $\delta>0$ be small number (how small
    $\delta$ is will be determined later) and let $B_\delta(X)$ be the
    ball in $\Ss^2$ of radius $\delta$ and center at $X$.  Fixing
    $\delta$, there exists a positive number $\e(\delta)>0$ such that
    $A_{\e}\big(\Ss^{2} \setminus B_\delta(X)\big) \subset
    B_\delta(-M_\e)$,
    $\norm{A_{\e}t_{\e}A_{\e}^{-1} - S_{\alpha}}_{C^{0}} < 0.01\delta$
    and
    $A^{-1}_{\e}S_{\alpha}\big(B_{1.01\delta}(-M_\e)\big) \subset
    B_\delta(X)$, for any $\e\in\big(0,\e(\delta)\big)$.

    And so we have

    \begin{displaymath}
      \begin{split}
        t_{\e} \big(\Ss^{2} \setminus B_\delta(X)\big) & =
        A^{-1}_{\e}\circ A_{\e}t_{\e}A_{\e}^{-1}\circ
        A_{\e}\big(\Ss^{2} \setminus B_\delta(X)\big) \\
        &\subset A^{-1}_{\e} \circ
        A_{\e}t_{\e}A_{\e}^{-1}\big(B_\delta(-M_\e)\big) \\
        & \subset A^{-1}_{\e} \circ
        S_{\alpha}\big(B_{1.01\delta}(-M_\e)\big)\subset B_\delta(X).
      \end{split}
    \end{displaymath}
    
    Now, for $\e>0$ sufficiently small we have $s_{\e}$ is close
    enough to the rotation $R_{\theta}\SO(3)$ constructed in the proof
    of Lemma \ref{confor}, so for such $\e$ we can also assume that
    $s_{\e}\big(B_{\delta}(X)\big)\subset\Ss^2 \setminus B_\delta(X)$
    and thus we get 
    \begin{displaymath}
      s_{\e}t_{\e}\big(\Ss^{2} \setminus B_\delta(X)\big) \subset \Ss^{2}
      \setminus B_\delta(X),
    \end{displaymath}
    where this last inclusion is strict and therefore $s_{\e}t_{\e}$
    has infinite order.
  \end{proof}
\end{proof}

Now we can finish the proof of Theorem~\ref{subexcase}:

\begin{proof}[Proof of Theorem~\ref{subexcase}]
  Let $s,t$ be two arbitrary elements of the generating set $S$. We
  can suppose that $s,t$ do not have a common fixed point, because
  otherwise they generate a finite cyclic group and in such a case we
  could reduce the number of generators.

  By the proof of Proposition \ref{2gen}, the conjugates
  $s_{\e}, t_{\e}$ converge to Möbious transformations
  $A_s, A_t\in\PSL_2(\C)$, as $\e\to 0$, generating a finite group. It
  also follows from the proof of Proposition \ref{2gen} that $A_s$ and
  $A_t$ cannot have a common fixed point.

  We will show that there exists a sub-sequence $\e_n \to 0$ such that
  for any $h \in S$, the sequence of conjugates
  $h_{\e_n} := g_{\e_n}^{-1}hg_{\e_n}$ converges to a non-trivial
  Möbious map $A_h\in\PSL_2(\C)$. To show this, let $p_{\e},q_{\e}$
  denote the fixed points of $h_{\e}$. Arguing as in Lemma
  \ref{confor}, if for some $h \in S$ there is a sub-sequence
  $h_{\e_n}$ such that $p_{\e_n}$ and $q_{\e_n}$ converge to two
  different points $p,q$, then $h_{\e_n}$ must converge to a finite
  order element $A_h\in\PSL_2(\C)$ fixing $p$ and $q$.

  If no such sub-sequence of $h_{\e}$ exists, then the sequences of
  points $p_{\e}, q_{\e}$ have a common limit point $X\in\Ss^2$. By
  possibly replacing of $A_s$ with $A_t$, one can suppose that $X$ is
  not a fixed point of $A_s$. Therefore, we are in the same situation
  as in Proposition \ref{hard} for the elements $s_{\e},h_{\e}$.
  Applying the very same argument we obtain an element of $G$ of
  infinite order, getting a contradiction.

  In conclusion, via a diagonal argument, we can find a sub-sequence
  $\e_n \to 0$ so that each of the conjugates $h_{\e_n}$ converges to
  an element of $\PSL_2(\C)$. This implies that the conjugate
  $w_{\e_n}$ of every element $w \in G$ converges to an element $A_w$
  of $\PSL_2(\C)$. Moreover, since the orders of $w$ and $A_w$
  coincide, by Lemma \ref{cantat1} we know that $A_w\neq id$ provided
  $w \neq id$.

  The previous discussion implies there is an injective homomorphism
  $\Phi\colon G\to\PSL_2(\C)$ sending $g$ to $A_g$.  As $G$ is
  periodic and $\PSL_2(\C)$ is linear, the group $G$ must be finite.
\end{proof}

\section{The exponential growth case}
\label{sec:exponential-growth-case}

In this section we finish the proof of Theorem~\ref{thm:main}. By
Theorem \ref{subexcase}, we can now assume that $G$ has exponential
growth.

The idea of the proof under this additional hypothesis goes as
follows: as $G$ does not have sub-exponential growth, there is a
constant $c\in (0,1)$ and a sequence of natural numbers
$n_j\to\infty$, as $j\to\infty$, such that given any point
$x\in\Ss^2$, by the pigeonhole principle for every $j\in\N$, there are
two elements $g_j, h_j\in G$, with $\abs{g_j}_S,\abs{h_j}_S\leq n_j$
and such that $d\big(g_j(x),h_j(x)\big)<c^{-n_j}$.

On the other hand, by Lemma \ref{ideal2} for any $\e>0$ and $j$
sufficiently large, the derivatives of $g_j$ and $h_j$ are bounded by
$e^{\e {n_j}}$. Therefore, the element $f_j := h_j^{-1}g_j$ move $x$
exponentially close to itself. Then, since the group $G$ has bounded
exponent, the orbit $\{f_j^i(x)\}_{i=1}^{k}$ of $x$ has exponentially
small diameter, where $k\in\N$ is a exponent for the whole group.

So, we will prove that this implies $f_j$ has a fixed point
exponentially close to $x$. This argument also apply if instead of a
single point $x$ we consider a finite collection of points
$x_i\in\Ss^2$. This implies there are non-trivial elements of $G$ with
as many fixed points as one wants, contradicting the fact that every
element of $G$ is conjugate to some rotation.

\begin{lemma}
  \label{excase}
  Let us suppose $G<\Diff\infty{}(\Ss^2)$ is a finitely generated
  periodic group of bounded exponent and has sub-exponential growth of
  derivative, \ie the conclusion of Lemma \ref{ideal2} holds.

  Then, $G$ does has sub-exponential growth.
\end{lemma}

\begin{remark}
  It is important to notice that this is the only part of the of
  Theorem \ref{thm:main} where we use the bounded exponent assumption
  on the group $G$.
\end{remark}

Observe that Theorem \ref{thm:main} just follows as a straightforward
combination of Theorem \ref{subexcase} and Lemma \ref{excase}. We will
now begin the proof of Lemma \ref{excase}. Throughout the proof we
will use the classical Vinogradov's ``O'' notation which states that
given two sequences $\{f_n\},\{g_n\}$ we have $f_n = O(g_n)$ if there
exists a constant $C, n_0>0$ such that $\frac{1}{C}f_n<g_n<Cf_n$ for
$n \geq n_0$.
\begin{proof}[Proof of Lemma~\ref{excase}]
  Reasoning by contradiction, let us suppose $G$ has exponential
  growth.

  Let $S$ be a finite generating set of $G$. By replacing $S$ with a
  larger generating set, we can suppose there is a sequence of natural
  numbers $\{n_j\}_j$, with $n_j\to\infty$ and such that
  \begin{displaymath}
    \sharp\big\{g\in G : \abs{g}\leq n_j\big\} \geq 2^{7n_j},
    \quad\forall j\geq 1.
  \end{displaymath}

  Observe the group $G$ naturally acts on
  $(\Ss^2)^3 := \Ss^2\times\Ss^2\times\Ss^2$ via the diagonal
  action. Then consider a fixed triple of distinct points of $\Ss^2$
  $\bar{x} =(x_1,x_2,x_3) \in (\Ss^2)^3$, let us define
  \begin{displaymath}
    O_j := \left\{\big(g(x_1),g(x_2),g(x_3)\big) \in (\Ss^2)^3 :  g\in
      G,\ \abs{g}_S\leq n_j\right\},
  \end{displaymath}
  and observe the $G$-orbit of $\bar{x}$ is equal to $\bigcup_j O_j$.

  Notice that if we endow $(\Ss^2)^3$ with the product Riemannian
  structure, the volume of a ball of radius $2^{-n}$ in $(\Ss^2)^3$ is
  $O(2^{-6n})$.

  Then, as there are at least $2^{7n_j}$ elements in
  $B_{n_j}:=\{g\in G :\abs{g}_S\leq n_j\}$, by the pigeonhole
  principle, for each $j\geq 1$ there are at least two different
  elements $g_j,h_j$ in $B_{n_j}$ such that their corresponding images
  of $\bar{x}\in(\Ss^2)^3$ satisfies $g_j(\bar x), h_j(\bar x)\in O_j$
  and $d\big(g_j(\bar x), h_j(\bar x)\big) < 2^{-n_j}$, for every
  $j\geq 1$.

  On the other hand, by Lemma \ref{ideal2} we know that given any
  $\e>0$, there exists $N_2$ so for any $n \geq N_2$ and any element
  $g \in B_n$, it holds $\norm{D(g)}\leq e^{\frac{\e}{2} n}$. If we
  define $f_j := h_j^{-1}g_j$, we have $f_j\in B_{2n_j}$ and so
  $\norm{D(f_j)} \leq e^{\e n_j}$, for every $n_j\geq N_2$.  So, for
  each $i\in\{1,2,3\}$ and any $j\geq 1$ such that $n_j\geq N_2$, we
  have
  \begin{equation}
    \label{eq:fj-small-displace-xi}
    d\big(f_j(x_i), x_i\big) \leq \norm{D(h_j^{-1})} d\big(g_j(\bar{x}), 
    h_j(\bar{x})\big) \leq e^{\e n_j}2^{-n_j}.
  \end{equation}

  Then we will show that \eqref{eq:fj-small-displace-xi} implies that,
  for $n_j$ sufficiently large, $f_j$ has at least three fixed points,
  and consequently it is the identity, contradicting the fact that
  $g_j$ and $h_j$ were different.
  
  \begin{proposition}
    \label{fixedp}
    For each $j\geq 1$ with $n_j\geq N_2$ and each $i\in\{1,2,3\}$,
    there exists $p_{i,j}\in\Ss^2$ such that $f_j(p_{i,j})=p_{i,j}$
    and 
    \begin{displaymath}
      d(x_{i},p_{i,n}) = O(2^{-n}e^{(k+2)\e n}). 
    \end{displaymath}

    Consequently, $f_j=id$, for $j$ sufficiently large.
  \end{proposition}

  \begin{proof}[Proof of Proposition~\ref{fixedp}]
    It is clearly enough to prove the statement just for the point
    $x_1$.  Let $S_j$ be the shortest geodesic segment on $\Ss^2$
    joining the points $f_j(x_1)$ and $x_1$. Then we define
    \begin{displaymath}
      K_j := \bigcup_{i=1}^k f_j ^i(S_j),
    \end{displaymath}
    where $k$ is the exponent of $G$, \ie $g^k=id$, for every
    $g\in G$.  Observe that $K_j$ is a compact $f_j$-invariant
    set. Also, as $\norm{D(f_j)} \leq e^{\e n_j}$, the compact set
    $K_j$ has diameter of order $ O\big(2^{-n_j}e^{\e (k+1)n_j}\big)$.

    Therefore, for each $j$ there exists an open disk
    $D_j\subset\Ss^2$ whose radius is of order
    $O(2^{-n}e^{\e (k+1)n})$ and containing $K_j$. If $n_j$ is
    sufficiently large, then the set $\Ss^2 \setminus K_j$ has exactly
    one connected component $R_j$ containing the set
    $\Ss^2 \setminus D_j$. Therefore, since
    $\norm{D(f_j)} \leq e^{\e n_j}$, by area considerations, the
    connected component $R_j$ must be $f_j$-invariant, provided $n_j$
    is large enough.

    So, the set $K_j':= S^2 \setminus R_j$ is an $f_j$-invariant
    non-separating continuum and by Cartwright-Littlewood Theorem (see
    for instance \cite{BrownShortProof}), there exists a fixed
    $p_{1,j}\in K'_j$ of $f_j$.

    As $K_j'\subset D_j$, we have that the distance between $p_{1,j}$
    and $x_1$ should be of order $O(2^{-n}e^{\e(k+1)n})$.
  \end{proof}
\end{proof}

\section{The Burnside problem on $\T^q$ and hyperbolic manifolds}
\label{sec:burnside-problem-Td-hyperbolic-mfds}

In this section we prove
Theorems~\ref{thm:Burnside-problem-hyperbolic-mfd} and
\ref{thm:Burnside-problem-Td}, which extends previous results of
Guelman and Liousse~\cite{GuelmanLiousseBurnsideMeasurePres,
  GuelmanLiousseBurnsideGroupsHomeo} to higher dimensions.

\subsection{Displacement functions and Newman's theorem}
\label{sec:displ-funct-newman-thm}

In this paragraph we introduce some notation we shall need in the
proofs of Theorems \ref{thm:Burnside-problem-hyperbolic-mfd} and
\ref{thm:Burnside-problem-Td}.

Let $(X,d)$ be an arbitrary metric space. We define the
\emph{displacement function} $\Dis\colon\Homeo{}(X)\to [0,\infty]$ by
\begin{equation}
  \label{eq:Disp-function-def}
  \Dis(f):=\sup_{x\in X} d\big(f(x),x\big), \quad\forall f\in\Homeo{}(X),
\end{equation}
where $\Homeo{}(X)$ denotes the group of homeomorphisms of $X$.

Then, let us recall a classical result due to
Newman~\cite{NewmanThmPerTrans}:

\begin{theorem}
  \label{thm:newman}
  Let $M$ be a connected complete manifold, $d$ be a distance function
  compatible with the topology of $M$ and $k\geq 1$ be a natural
  number. Then, there exists a real number $\varepsilon(M,d,k)>0$,
  depending on $M$, $d$ and $k$, such that for every $f\in\Homeo{}(M)$
  satisfying $f^k=id$ and $f\neq id$, it holds $\Dis(f)\geq\varepsilon$.
\end{theorem}

An easy consequence of Theorem~\ref{thm:newman} is the following
\begin{corollary}
  \label{cor:newman}
  Let $d$ denote the standard Euclidean distance on $\R^q$ and
  $f\in\Homeo{}(\R^q)$ be a homeomorphism such that $f\neq id$ and
  $f^k=id$, for some $k\geq 1$. Then it holds $\Dis(f)=\infty$.
\end{corollary}

\begin{proof}[Proof of Corollary~\ref{cor:newman}]
  Let $\varepsilon:=\varepsilon(\R^q,d,k)>0$ be the real constant
  given by Theorem \ref{thm:newman} for $M=\R^q$, and let us suppose
  $\Dis(f)<\infty$. Now, if we define $f_\lambda\in\Homeo{}(\R^q)$ by
  \begin{displaymath}
    f_\lambda(x):=\lambda f\big(\lambda^{-1} x\big),\quad\forall
    x\in\R^q,
  \end{displaymath}
  with $\lambda:=\varepsilon (2\Dis(f))^{-1}$, then it holds
  \begin{displaymath}
    \Dis(f_\lambda)=\sup_{x\in\R^q} \norm{f_\lambda(x)-x} =
    \sup_{x\in\R^q} \lambda \norm{f\big(\lambda^{-1} x\big)
      -\lambda^{-1}x} = \lambda \Dis(f) = \frac{\varepsilon}{2},
  \end{displaymath}
  and $f_\lambda^k=id$. Then, by Theorem~\ref{thm:newman}
  it holds $f_\lambda=id$, contradicting the fact that $f\neq id$. 
\end{proof}

\subsection{Proof of Theorem \ref{thm:Burnside-problem-hyperbolic-mfd}}
\label{sec:proof-theorem-Burnside-prb-hyp-mfd}

Let $M$ be compact hyperbolic manifold, $d$ denote the distance
function induced by its Riemannian structure and let $G<\Homeo{}(M)$
be a finitely generated periodic group.

By Mostow rigidity theorem \cite{MostowQuasiconformalMappings} we know
that each homeomorphism of $M$ is homotopic to an isometry and
$\Isom(M)$, the group of isometries of $M$, is a finite group (see for
instance \cite[Theorem C.5.6]{BenedettiPetronioLecturesHypGeo} for
details).

So, if we write $\Homeo0(M)$ for the subgroup of homeomorphisms of $M$
which are homotopic to the identity and we define
$G_0:=G\cap\Homeo0(M)$, we conclude $G_0$ has finite index in $G$. The
rest of the proof consists in showing that $G_0$ is a singleton,
containing just the identity. 

To do that, let $\pi\colon\bb{H}^q\to M$ be the universal covering of
$M$ and write $\Deck(\pi)<\Isom(\bb{H}^q)$ for the group of
automorphisms of the covering. We shall write $\widetilde d$ for the
classical hyperbolic distance on $\bb{H}^q$.

Since $M$ is compact, we know that the every element of
$\Deck(\pi)\setminus\{id\}$ has trivial centralizer (\ie the
centralizer is the infinite cyclic group). This implies that for every
$f\in\Homeo0(M)$, there exists a unique lift $\tilde
f\in\Homeo{}(\bb{H}^q)$ such that
\begin{displaymath}
  \tilde f\circ\gamma = \gamma\circ \tilde f,
  \quad\forall\gamma\in\Deck(\pi). 
\end{displaymath}
Observe that this lift $\tilde f$ can be also characterized as the
final element of an identity homotopy of $\bb{H}^q$ (\ie a homotopy
starting at the identity map) that is a lift of an identity homotopy
on $M$ for $f$. In particular, if
$D\colon\Homeo{}(\bb{H}^q)\to [0,\infty]$ denotes the displacement
function given by \eqref{eq:Disp-function-def}, then it clearly holds
\begin{equation}
  \label{eq:D-f-finite}
  \Dis(\tilde f)<\infty, \quad\forall f\in\Homeo0(M).
\end{equation}
On the other hand, it is well known that
\begin{equation}
  \label{eq:D-gamma-infinite}
  \Dis(\gamma)=\infty,
  \quad\forall\gamma\in\Isom(\bb{H}^q)\setminus\{id\}. 
\end{equation}

Then, let us define
\begin{displaymath}
  \widetilde G_0:=\left\{\tilde g : g\in G\right\}. 
\end{displaymath}
Notice that, since $\widetilde {f\circ g} = \tilde f\circ \tilde g$,
for every $f,g\in\Homeo0(M)$, we have that $\widetilde G_0$ is \emph{à
  priori} just a semi-group and not a group. But we shall show that
$\widetilde G_0=\{id\}$.

In fact, since $G$ is periodic, for every $g\in G$, there exist
$\gamma_g\in\Deck(\pi)$ and $n_q\in\N$ such that
$\tilde g^{n_g}=\gamma_g$. So, by \eqref{eq:D-f-finite} we know that
every $g\in\widetilde G$ it holds
\begin{displaymath}
  \Dis(\gamma_g)=D\big(\tilde g^{n_g}\big) \leq n_g\Dis(\tilde g) <\infty.
\end{displaymath}
And then, by \eqref{eq:D-gamma-infinite} we conclude that
$\gamma_g=id$, for every $g\in G_0$. So, $\widetilde G_0$ is a
periodic group as well.

Finally, combining Corollary~\ref{cor:newman} and
\eqref{eq:D-f-finite} (and observing that $\bb{H}^q$ is homeomorphic
to $\R^q$), we conclude that $\widetilde G_0=\{id\}$. So, $G_0$ is a
singleton as well and hence, $G$ is finite.

\subsection{Proof of Theorem~\ref{thm:Burnside-problem-Td}}
\label{sec:proof-theorem-Burnside-Td}

Let $\pi\colon\R^q\to\T^q:=\R^q/\Z^q$ be the natural projection
map. For each $f\in\Homeo{}(\T^q)$, there is a unique
$A_f\in\GL_q(\Z)$ such that the map
$\Delta_{\tilde f}:=\tilde f - A_f\colon\R^q\to\R^q$ is
$\Z^q$-periodic, for any lift $\tilde f\colon\R^q\carr$ of $f$.  The
map
\begin{displaymath}
  \mathcal{I} : \Homeo{}(\T^q)\ni f\mapsto A_f\in\GL_q(\R)
\end{displaymath}
is clearly a group homomorphism. Then we define
$\Homeo0(\T^q):=\ker\mathcal I$.

Let $G<\Homeo{\mu}(\T^q)$ denote a finitely generated periodic
subgroup. Then, $\mathcal{I}(G)<\GL_q(\R)$ is a finitely generated
periodic group. By Schur's theorem we know $\mathcal{I}(G)$ is
finite. So, the sub-group
\begin{displaymath}
  G_0:=G\cap \Homeo0(\T^q)<G
\end{displaymath}
has finite index in $G$. Thus, it is enough to show $G_0$ is finite.

Then, let us consider the group
$\Homeo{0,\mu}(\T^q):=\Homeo\mu(\T^q)\cap\Homeo0(\T^q)$ and the
\emph{$\mu$-mean rotation vector}
$\rho_\mu\colon\Homeo{0,\mu}(\T^q)\to\T^q$ given by
\begin{displaymath}
  \rho_\mu(f):=\pi\left(\int_{\T^q} \Delta_{\tilde f}\dd\mu\right),
  \quad\forall f\in\Homeo{0,\mu}(\T^q),
\end{displaymath}
where $\tilde f\colon\R^q\carr$ is any lift of $f$ and the function
$\Delta_{\tilde f}=\tilde f-id_{\R^q}$ is considered as an element of
$C^0(\T^q,\R^q)$. Notice that $\rho_\mu$ is well-defined, does not
depend on the choice of lift.

We claim $\rho_\mu$ is a group homomorphism. To prove this, let
$f,g\in\Homeo{0,\mu}(\T^q)$ be two arbitrary homeomorphisms,
$\tilde f,\tilde g\colon\R^q\carr $ two lifts of $f$ and $g$,
respectively. Then we have:

\begin{displaymath}
  \begin{split}
    \rho_\mu(f\circ g) &= \pi\left(\int_{\T^q}\Delta_{\tilde
        f\circ\tilde g}\dd\mu \right) = \pi\left(\int_{\T^q}
      \Delta_{\tilde f}\circ g + \Delta_{\tilde g}\dd\mu\right) \\
    &= \pi\left(\int_{\T^q}\Delta_{\tilde f}+\Delta_{\tilde g}\dd\mu
    \right) = \rho_\mu(f) + \rho_\mu(g),
  \end{split}
\end{displaymath}
where the first equality is consequence of the identity
$\Delta_{\tilde f\circ\tilde g}=\Delta_{\tilde f}\circ g +
\Delta_{\tilde g}$ and the second one follows from the fact that $\mu$
is $g$-invariant.

So, $\rho_\mu(G_0)<\T^q$ is a finitely generated periodic group; and
since $\T^q$ is abelian, $\rho_\mu(G_0)$ is finite. This implies
$G_{0,0}:=G_0\cap\ker\rho_\mu$ is a finite index sub-group of
$G_0$. Thus, it is enough to show $G_{0,0}$ is finite, and we shall
indeed show it just reduces to the identity map.

In order to do that, let $f$ be an arbitrary element of
$G_{0,0}$. Notice there is a unique lift $\tilde f\colon\R^q\carr$
such that
\begin{displaymath}
  \int_{\T^q}\Delta_{\tilde f}\dd\mu = 0.
\end{displaymath}
Since $G_{0,0}$ is periodic, there exists $n\geq 1$ such that
$f^n=id$, and consequently, there is $\boldsymbol{p}\in\Z^q$ such that
$\tilde f^n(x)=x+\boldsymbol{p}$, for every $x\in\R^q$. However, since
$\mu$ is $f$-invariant, it holds
\begin{displaymath}
  \boldsymbol{p}=\int_{\T^q}\Delta_{\tilde f^n}\dd\mu = \int_{\T^q}
  \sum_{j=0}^{n-1}\Delta_{\tilde f}\circ f^j\dd\mu =
  n\int_{\T^q}\Delta_{\tilde f}\dd\mu = 0.
\end{displaymath}
So, $\tilde f$ is a periodic homeomorphism of $\R^q$ and since
$\Delta_{\tilde f}$ is $\Z^q$-periodic, it holds
\begin{displaymath}
  \Dis(\tilde f)=\sup_{x\in\R^q} \norm{\Delta_{\tilde f}(x)} <\infty,
\end{displaymath}
where $\Dis(\tilde f)$ is the displacement function given by
\eqref{eq:D-f-finite}. Then, by Corollary~\ref{cor:newman} we get
$\tilde f$ equals the identity Theorem~\ref{thm:Burnside-problem-Td}
is proven.

\bibliographystyle{amsalpha} \bibliography{references}

\end{document}